\documentclass[11pt,reqno]{amsart}
\usepackage[a4paper,margin=30mm]{geometry}
\usepackage{amsmath,amssymb,mathtools}
\usepackage{microtype,booktabs}
\usepackage[colorlinks=true,linkcolor=blue,citecolor=blue,urlcolor=blue]{hyperref}
\hypersetup{bookmarksdepth=1}

\newtheorem{theorem}{Theorem}[section]
\newtheorem{proposition}[theorem]{Proposition}
\newtheorem{lemma}[theorem]{Lemma}

\theoremstyle{definition}

\theoremstyle{remark}
\newtheorem{remark}[theorem]{Remark}

\numberwithin{equation}{section}

\newcommand{\R}{\mathbb{R}}

\newcommand{\Q}{\mathbb{Q}}
\newcommand{\Z}{\mathbb{Z}}
\newcommand{\PP}{\mathbb{P}}
\newcommand{\calA}{\mathcal{A}}
\newcommand{\cI}{\mathrm{I}}
\newcommand{\cII}{\mathrm{II}}
\newcommand{\cIII}{\mathrm{III}}
\newcommand{\cIV}{\mathrm{IV}}
\DeclareMathOperator{\rk}{rank}
\DeclareMathOperator{\NS}{NS}

\DeclareMathOperator{\Aut}{Aut}

\title{A Linear Bound for the Second Betti Number of Hyperk\"ahler Manifolds}
\author{Jie Fu}
\address{Jie Fu, Qiuzhen College, Tsinghua University, Beijing, China}
\email{fu-j21@mails.tsinghua.edu.cn}
\author{Shihao Wang}
\address{Shihao Wang, Qiuzhen College, Tsinghua University, Beijing, China}
\email{wangshihao25@mails.tsinghua.edu.cn}
\author{Zhiwei Zheng}
\address{Zhiwei Zheng, Yau Mathematical Sciences Center, Tsinghua University, Beijing, China}
\email{zhengzhiwei@mail.tsinghua.edu.cn}
\date{}

\begin{document}

\begin{abstract}
Assuming the correctness of Proposition~7.6 of OpenAI's \emph{Projective-space bases of Lagrangian fibrations}~\cite{base}, we establish the linear bound $b_2(X)\le 8n+13$ for every smooth compact irreducible hyperk\"ahler manifold $X$ of complex dimension $2n\ge 10$. We restrict a Lagrangian fibration to a general line in its base and bound the ramification of its first cohomology local system from below using transcendental cohomology and from above using intersection theory on a compactified moduli space of abelian varieties.
\end{abstract}

\maketitle


\section{Introduction}\label{sec:introduction}

A \emph{smooth compact irreducible hyperk\"ahler manifold} (an IHS manifold for short) $X$ of complex dimension $2n$ is a simply connected compact K\"ahler manifold such that $H^0(X,\Omega_X^2)$ is spanned by an everywhere nondegenerate holomorphic two-form. We write $b_j(X) = \dim_\Q H^j(X,\Q)$ for the $j$-th Betti number of $X$.

Beauville~\cite{Beauville} constructed two series of IHS manifolds. Namely, for each $n\ge 2$, there exists two deformation types in complex dimension $2n$, represented by the Hilbert scheme $\mathrm{K3}^{[n]}$ of $n$ points on a K3 surface, with $b_2 = 23$, and the generalized Kummer variety $\mathrm{Kum}_n$, with $b_2 = 7$. Beyond these two series, O'Grady discovered two sporadic examples in complex dimensions six and ten, denoted $\mathrm{OG6}$ and $\mathrm{OG10}$ whose second Betti numbers $8$ and $24$, respectively.

A central question in hyperk\"ahler geometry is to bound $b_2(X)$ from above in terms of the dimension $2n$. Every IHS surface is a K3 surface with $b_2 = 22$. For fourfolds, Guan~\cite[Main Theorem]{Guan} proved the optimal inequality $b_2 \le 23$. For sixfolds ($2n=6$), Sawon~\cite{saw-six} deduced $b_2 \le 23$ conditionally on a hypothesis concerning the Looijenga--Lunts--Verbitsky (LLV) decomposition of cohomology, and the unconditional sharp bound $b_2 \le 23$ was established in~\cite{FWWZ}.

In general dimension, two conditional programs have been proposed to bound $b_2(X)$. Green, Kim, Laza, and Robles~\cite{GKLR} computed the LLV decomposition for all known deformation types and conjectured restrictions on the highest weights that may appear; Kim and Laza~\cite[Main Theorem]{KL} showed that this conjecture yields a dimension-dependent upper bound on $b_2$. From the topological side, Sawon~\cite{saw-top} proposed a sign conjecture for Rozansky--Witten invariants and demonstrated that it likewise implies upper bounds on $b_2$ in all dimensions.

The purpose of this paper is to prove an upper bound on $b_2(X)$ that grows linearly with $n$.

\begin{theorem}\label{thm:main}
Let $X$ be an IHS manifold of complex dimension $2n$ with $n \ge 5$. Then
\begin{equation*}
  3 \le b_2(X) \le 8n + 13.
\end{equation*}
\end{theorem}

The proof proceeds as follows. Assuming $b_2(X) \ge 5$, we deform $X$ to a projective IHS manifold of Picard number two carrying a Lagrangian fibration. We use OpenAI~\cite[Proposition~7.6]{base} only to conclude that the base has at worst quotient singularities. We apply~\cite[Theorem~1.2]{liu} to identify the base with $\PP^n$.

Let $S$ be the set of singular values on a general line $\ell\subset\PP^n$. Let $\mathbb V$ be the local system of first rational cohomology groups of the smooth fibres over $\ell\setminus S$, with local monodromy $M_s$ around each $s\in S$. The \emph{ramification} of $\mathbb V$ along $\ell$ is $\kappa=\sum_{s\in S}\rk(M_s-\operatorname{Id})$. It is the sum of the codimensions of the fixed subspaces of the $M_s$. We prove that
\begin{equation}\label{eq:intro-ramification}
  b_2(X) + 4n - 2 \le \kappa \le 12n + 11.
\end{equation}
For the lower bound, the transcendental part of the projective deformation's second cohomology has dimension $b_2(X)-2$ and embeds into the intersection cohomology group $IH^1(\ell,\mathbb V)$, whose dimension is $\kappa-4n$.

For the upper bound, we apply Zarhin's quaternion trick to produce an associated family of principally polarized abelian varieties, compute its normalized Hodge degree and its boundary multiplicities in the perfect cone compactification, and verify that the resulting period curve is not contained in the exceptional locus of $12\lambda - D$. Invoking Shepherd-Barron's positivity theorems~\cite{sb1,sb2} then yields $\kappa \le 12n + 11$. Combining both sides of~\eqref{eq:intro-ramification} immediately gives $b_2(X) \le 8n + 13$.

\begin{remark}
For IHS eightfolds ($n=4$), the same argument gives $\kappa\le12n+12$ and further
\[
  3\le b_2(X)\le8n+14=46.
\]
The assumption $n\ge5$ is used to prove the strict positivity of $\deg_L p^*(12\lambda-D)$. In dimension eight, nefness still gives $\deg_L p^*(12\lambda-D)\ge0$, so the upper bound for $\kappa$ increases by one. 

We will prove the sharper bound $b_2(X)\le24$ for IHS eightfolds in a subsequent paper. Together with previous results, this gives an explicit upper bound for the second Betti number of IHS manifolds in every dimension.
\end{remark}

\medskip
\noindent\textbf{Acknowledgements.}
We thank Bin Wang for related discussions. As we completed the first draft, we learned that at least one other group had also used AI to obtain a similar linear bound on the second Betti number. Our work and theirs were carried out independently of one another.

\medskip
\noindent\textbf{Use of AI tools.}
We developed the main arguments of this paper over multiple rounds of discussion with GPT-6 Astra and GPT-6.2 Astra. We independently checked and finalized the entire manuscript and take full responsibility for its contents.

\section{Reduction to a projective Lagrangian fibration}\label{sec:reduction}

\subsection{Deformation to a Lagrangian fibration of Picard rank two}

Let $X$ be an IHS manifold of complex dimension $2n$, with $n \ge 5$. Since Theorem~\ref{thm:main} holds trivially when $b_2(X) \le 4$, we assume throughout the proof that $b_2(X) \ge 5$.

Let $q$ denote the primitive integral Beauville--Bogomolov--Fujiki quadratic form on $H^2(X,\Z)$. It has signature $(3, b_2(X)-3)$ and satisfies the Fujiki relation
\begin{equation}\label{eq:fujiki-general}
  \int_X \gamma^{2n} = c_X \, q(\gamma)^n, \qquad \gamma \in H^2(X,\R),
\end{equation}
for a positive constant $c_X > 0$ depending only on the deformation class of $X$~\cite[Th\'eor\`eme~5]{Beauville}.

\begin{proposition}\label{prop:model}
Let $X$ be an IHS manifold of complex dimension $2n$ with $b_2(X) \ge 5$. Then $X$ is deformation equivalent to an IHS manifold $X'$ of Picard number $2$ admitting a Lagrangian fibration $f \colon X' \to \PP^n$.
\end{proposition}
\begin{proof}
Since $b_2(X) \ge 5$, the integral quadratic form $q$ of signature $(3, b_2(X)-3)$ is indefinite in at least five variables. By Meyer's theorem, $q$ has zero in a primitive nonzero integral cohomology class $\eta \in H^2(X,\Z)$. We deform $X$ to a very general period in $\eta^\perp$, so that $\eta$ is of type $(1,1)$ and generates the Picard group. Since there are no negative-square classes in the Picard group, $\eta$ is nef after a choice of sign. By Engel--Mauri~\cite[Theorem~1.1]{em}, Matsushita's theorem~\cite{Mat01}
and its K\"ahler version~\cite[Theorem~2.10]{KamenovaLehn},
a positive multiple of $\eta$ defines a Lagrangian fibration.

By~\cite[Corollaries~A.3--A.4]{efgms}, this Lagrangian fibration admits a projective deformation. Let $H$ be an ample class on the projective total space, and let $\ell$
be the pullback of an ample class from the base. By~\cite[Corollary~A.4]{efgms}, the polarized fibration deforms locally over the locus where $H$ and $\ell$ remain of type $(1,1)$. For a very general member $X'$ of this family, $\NS(X')_\Q=\Q H\oplus\Q\ell$, so $X'$ has Picard number two. Finally, the projective-space base theorem~\cite[Theorem~1.1]{base} identifies the base of this fibration with $\PP^n$.
\end{proof}

Since smooth deformations preserve $b_2(X)$, Proposition~\ref{prop:model} allows us to replace $X$ by $X'$. We write $f \colon X \to \PP^n$ for this projective model of Picard number two in what follows.

\subsection{Singular fibres and local monodromy on a general line}

Let $U \subset \PP^n$ denote the smooth locus of $f$, and let $\Delta = \PP^n \setminus U$ be the discriminant locus. Choose a general projective line $L \subset \PP^n$, put $S = L \setminus U = L \cap \Delta$, and denote by
\begin{equation*}
  j \colon L^\circ = L \setminus S \hookrightarrow L
\end{equation*}
the open inclusion. Write $V = (R^1 f_* \Q)|_{L^\circ}$ for the weight-one rational variation of Hodge structure of rank $2n$ on $L^\circ$, and let $T_s \in \Aut(V_{t_0})$ denote the local monodromy operator at $s \in S$.

A general line $L\subset\PP^n$ meets each divisorial component of $\Delta$ transversely at general points and avoids the components of $\Delta$ of codimension at least two. By the theorem of Kim, Oguiso, and Shinder~\cite[Corollary~1.2]{kos}, the singular fibres of $f$ over $S$ are reduced and have type $\cI, \cII, \cIII, \cIV$. By Kim's classification~\cite[Theorem~5.1.1 and \S 6.2]{kim}, the local monodromy $T_s$ at $s \in S$ behaves as follows:
\begin{enumerate}
\item[\textup{(i)}] At a point of type $\cI$, the local monodromy $T_s$ is unipotent with a rank-one nilpotent logarithm, so $\rk(T_s - \operatorname{Id}) = 1$.
\item[\textup{(ii)}] At a point of type $\cII, \cIII$, or $\cIV$, the local monodromy $T_s$ is the identity outside a rank-two elliptic summand on which it has finite order $6, 4$, or $3$, respectively, with eigenvalues $e^{\pm 2\pi i / d}$ ($d \in \{6,4,3\}$). In each of these three additive types, $\rk(T_s - \operatorname{Id}) = 2$.
\end{enumerate}

Let $N_{\cI}, N_{\cII}, N_{\cIII}, N_{\cIV}$ denote the numbers of points in $S$ of types $\cI, \cII, \cIII, \cIV$, respectively. The local invariants of these four fibre types are summarized in Table~\ref{tab:fibres}.

\begin{table}[htbp]
\centering
\renewcommand{\arraystretch}{1.3}
\begin{tabular*}{\textwidth}{@{\extracolsep{\fill}}lcccc@{}}
\toprule
Fibre type & Monodromy action & Order $d_s$ & $\rk(T_s - \operatorname{Id})$ & Hodge degree drop\\
\midrule
$\cI$ & Unipotent (rank-$1$ logarithm) & $\infty$ & $1$ & $0$\\
$\cII$ & Elliptic finite order & $6$ & $2$ & $1/6$\\
$\cIII$ & Elliptic finite order & $4$ & $2$ & $1/4$\\
$\cIV$ & Elliptic finite order & $3$ & $2$ & $1/3$\\
\bottomrule
\end{tabular*}
\caption{Local monodromy invariants of singular fibres along a general line $L \subset \PP^n$.}
\label{tab:fibres}
\end{table}

We define the \emph{ramification} of $V$ along $L$ by
\begin{equation}\label{eq:ramification-def}
  \kappa \coloneqq \sum_{s \in S} \rk(T_s - \operatorname{Id}) = N_{\cI} + 2\bigl(N_{\cII} + N_{\cIII} + N_{\cIV}\bigr).
\end{equation}

\section{The Cohomological Lower Bound}\label{sec:lower}

In this section we prove the lower bound for the ramification $\kappa$ in~\eqref{eq:intro-ramification}.

Put $Y = f^{-1}(L)$ and $Y^\circ = f^{-1}(L^\circ)$. Let $\ell = f^* c_1(\mathcal{O}_{\PP^n}(1)) \in \NS(X)$ be the pullback of the hyperplane class from the base, and let $H \in \NS(X)$ be an ample divisor class on $X$. By Bertini's theorem, $Y$ is a smooth compact submanifold of $X$ of complex dimension $n+1$ with fundamental cohomology class $[Y] = \ell^{n-1} \in H^{2n-2}(X,\Q)$.

Consider the rational transcendental lattice
\begin{equation*}
  T \coloneqq \NS(X)_\Q^{\perp_q} \subset H^2(X,\Q).
\end{equation*}
Since $\rk \NS(X) = 2$ and $\NS(X)_\Q$ is spanned by the isotropic class $\ell$ and the ample class $H$ with $q(H,\ell) > 0$, the restriction of $q$ to $\NS(X)_\Q$ is nondegenerate. Consequently, $\dim_\Q T = b_2(X) - 2$, the bilinear form $q|_T$ is nondegenerate, and $T$ is a polarizable pure rational Hodge structure of weight two containing no nonzero rational $(1,1)$-classes.

\begin{lemma}\label{lem:transcendental-inject}
The restriction maps $T \to H^2(Y,\Q)$ and $T \to H^2(Y^\circ,\Q)$ are injective.
\end{lemma}
\begin{proof}
For any $\alpha, \beta \in T$, we have $q(\alpha,\ell) = q(\beta,\ell) = q(\alpha,H) = q(\beta,H) = 0$ and $q(\ell,\ell) = 0$. Polarizing the Fujiki relation~\eqref{eq:fujiki-general} on $X$ therefore gives
\begin{equation}\label{eq:fujiki-Y}
  \int_Y (\alpha|_Y) \smile (\beta|_Y) \smile (H|_Y)^{n-1} = \int_X \alpha \smile \beta \smile H^{n-1} \smile \ell^{n-1} = c \, q(\alpha,\beta) q(H,\ell)^{n-1},
\end{equation}
for a constant $c > 0$. Since $q(H,\ell) > 0$ and $q|_T$ is nondegenerate, relation~\eqref{eq:fujiki-Y} shows that $T$ injects into $H^2(Y,\Q)$.

By cohomological purity for the open submanifold $Y^\circ \subset Y$, the kernel of the restriction map $H^2(Y,\Q) \to H^2(Y^\circ,\Q)$ is spanned by the fundamental classes of irreducible divisorial components of the singular fibres $Y \setminus Y^\circ = f^{-1}(S)$. In particular, this kernel consists entirely of rational $(1,1)$-classes. Since the pure weight-two Hodge structure $T$ contains no nonzero rational $(1,1)$-classes, it intersects the kernel trivially and hence still injects into $H^2(Y^\circ,\Q)$.
\end{proof}

\begin{proposition}\label{prop:lower}
The ramification on a general line $L \subset \PP^n$ satisfies
\begin{equation*}
  \kappa \ge b_2(X) + 4n - 2.
\end{equation*}
\end{proposition}
\begin{proof}
Since $\dim_\Q T=b_2(X)-2$, it suffices to show that
\[
  \dim_\Q T \le h^1(L,j_*V)=\kappa-4n.
\]
We first compute $h^1(L,j_*V)$.

By the global invariant cycle theorem and $H^1(X,\Q) = 0$, the local system $R^1 f_* \Q$ has no nonzero global invariants over $U$. For a general line $L \subset \PP^n$, the Lefschetz theorem for hypersurface complements of Dimca and Papadima~\cite[Theorem~1]{dp} gives a surjection $\pi_1(L^\circ) \twoheadrightarrow \pi_1(U)$, so
\begin{equation*}
  H^0(L, j_* V) = H^0(L^\circ, V) = 0.
\end{equation*}
Intersection cohomology duality together with the polarization on $V$ likewise gives $H^2(L, j_* V) = 0$. At each singular point $s \in S$, the stalk $(j_* V)_s = \ker(T_s - \operatorname{Id})$ has dimension $2n - \rk(T_s - \operatorname{Id})$. Additivity of the Euler characteristic on $L = L^\circ \sqcup S$ therefore yields
\begin{equation*}
  \chi(L, j_* V) = 2n\bigl(2 - |S|\bigr) + \sum_{s \in S} \bigl(2n - \rk(T_s - \operatorname{Id})\bigr) = 4n - \kappa.
\end{equation*}
Since $H^0(L, j_* V) = H^2(L, j_* V) = 0$, it follows that $h^1(L,j_*V)=\kappa-4n$.

It remains to prove the dimension inequality. By Lemma~\ref{lem:transcendental-inject}, $T$ injects into $H^2(Y^\circ,\Q)$. We first show that its image restricts to zero
on every smooth fibre.

On each smooth fibre $X_t = f^{-1}(t)$ ($t \in L^\circ$), the holomorphic symplectic two-form vanishes because the fibres of $f$ are Lagrangian. Consequently, the image of $T$ in $H^2(X_t,\Q)$ has pure Hodge type $(1,1)$. Since polarizable pure rational Hodge structures are semisimple and $T$ has no nonzero rational $(1,1)$-classes, $T$ admits no nonzero quotient of type $(1,1)$. Thus $T$ restricts to zero on every smooth fibre over $L^\circ$.

By~\cite[Proposition~4.1]{Hwang} or~\cite[Theorem~1.7]{varmt}, the discriminant locus $\Delta \subset \PP^n$ possesses a divisorial component, so $S = L \cap \Delta$ is non-empty and the punctured line $L^\circ = L \setminus S$ is an affine curve. Because $L^\circ$ is affine of complex dimension one, $H^2(L^\circ, f_* \Q) = 0$. The Leray spectral sequence for $f|_{Y^{\circ}}$ therefore gives the following exact sequence of mixed Hodge structures
\[
0 \longrightarrow H^1(L^\circ,V)
  \longrightarrow H^2(Y^\circ,\Q)
  \longrightarrow H^0(L^\circ,R^2f_*\Q)
  \longrightarrow 0.
\]

By Lemma~\ref{lem:transcendental-inject} and $T$ being pure of weight $2$, $T$ embeds into the weight-two part $W_2 H^1(L^\circ, V)$. By the lowest-weight theorem of Peters and Saito and the description of the intersection complex on a smooth projective curve~\cite[Main Theorem~(b)]{ps}, $W_2 H^1(L^\circ, V)$ is the image of $H^1(L, j_* V)$.

Hence
\[
  b_2(X)-2=\dim_\Q T\le h^1(L,j_*V)=\kappa-4n,
\]
as required.
\end{proof}

\section{The upper bound for the ramification}\label{sec:upper}

\subsection{The Hodge degree}

Let $p: A\to L^\circ$ denote the relative Albanese abelian scheme of translations of the smooth fibres of $f$ over $L^\circ$, equipped with the polarization induced by the relatively ample class $H|_{Y^\circ}$. Since $Y^\circ$ is an $A$-torsor and translations act trivially on cohomology, there is a canonical isomorphism $R^1p_*\Q \simeq V$ of variations of Hodge structure.

After a finite branched base change $\pi \colon \widetilde{L} \to L$ from a smooth projective curve $\widetilde{L}$ such that the local monodromy of the pullback is unipotent at every point above $S$, $\widetilde{A}=A \times_{L^\circ} \pi^{-1}(L^\circ)$ extends to a semi-abelian scheme over $\widetilde{L}$ by Grothendieck's semistable reduction criterion
\cite[Exp.~IX, Prop.~3.2 and Cor.~3.8]{SGA7I}. Let $h$ denote the degree of the determinant of its semistable Hodge bundle $p_*\Omega^1_{\widetilde{A}/\widetilde{L}^\circ}$ on $\widetilde{L}$ divided by $\deg(\widetilde{L}/L)$. Because the unipotent Deligne extension commutes with finite ramified pullbacks and is invariant under isogeny~\cite[\S~1.5, (1.5.2)]{saito-schnell}, the rational number $h$ is independent of the choice of the cover $\widetilde{L}$ and is preserved under $L^\circ$-isogenies.

\begin{lemma}\label{lem:degree}
The normalized Hodge degree $h$ satisfies
\begin{equation}\label{eq:hodge-degree}
  h = n + 1 - \frac{N_{\cII}}{6} - \frac{N_{\cIII}}{4} - \frac{N_{\cIV}}{3}.
\end{equation}
\end{lemma}
\begin{proof}
Over the smooth locus $U$, contraction with the holomorphic symplectic two-form identifies the Hodge bundle of the fibration with the tangent bundle $T_{\PP^n}|_U$ of the base. On $L$, take the Deligne extension with residue eigenvalues in $(-1,0]$, and denote its extended Hodge bundle by $E$. By~\cite[\S~2, equation~(2.7) and the discussion after Lemma~2.6]{bs}, the above identification extends to an isomorphism $E \simeq T_{\PP^n}|_L$. Here relative Hard Lefschetz identifies the weight-$(2n-1)$ variation considered in~\cite{bs} with our weight-one variation up to a Tate twist. The determinant of this extended Hodge bundle therefore has degree
\begin{equation*}
  \deg_L(\det T_{\PP^n}|_L) = n + 1.
\end{equation*}

To compute $h$, we compare two Hodge line bundles on $\widetilde L$: the pullback of the Hodge line bundle on $L$ and the semistable Hodge line bundle on $\widetilde L$. They agree away from $\pi^{-1}(S)$,
so it suffices to compare their local generators at the points above $S$.

At a point $s$ of type $\cI$, the local monodromy is unipotent, so the Hodge bundle of the Deligne extension on $L$ pulls back to the Hodge bundle of the Deligne extension on $\widetilde L$ near every point above $s$. Thus there is no local correction to the degree. 

At a point $s \in S$ of type $\cII, \cIII$, or $\cIV$, let $t$ be a local coordinate on $L$ centered at $s$ and compute on the minimal local cover $t = u^d$, where $d = 6, 4, 3$ in the three respective types. On the good reduction model over the $u$-disk, let $g$ be the generator of the local Galois group given by $g(u)=\zeta_d u$, where $\zeta_d=e^{2\pi i/d}$. By Kim's local description for types $\cII$, $\cIII$, and $\cIV$~\cite[\S~6.2.1, Propositions~6.2.1--6.2.2]{kim}, we can choose a local frame $\omega,\omega_2,\ldots,\omega_n$
of the semistable Hodge bundle on the $u$-disk such that
\[
g^*\omega=\zeta_d\omega,
\qquad
g^*\omega_j=\omega_j \quad (2\le j\le n).
\] Since $g^*(u^{-1}\omega)=u^{-1}\omega$, the section $u^{-1}\omega$ descends over the punctured $t$-disk. The exponent $-1$ corresponds to the residue $-1/d\in(-1,0]$, so the pullback of the Hodge bundle of the $(-1,0]$-Deligne extension on $L$ has local frame $u^{-1}\omega,\omega_2,\ldots,\omega_n$.
The Hodge bundle on the good reduction model has local frame $\omega,\omega_2,\ldots,\omega_n$, whose determinant is $u$ times the determinant of the preceding frame. Thus the local degree correction is $-1$, giving a normalized correction of $-1/d$ at $s$. Any further ramified pullback multiplies the discrepancy by its ramification index, so the normalized drop is the same on the chosen cover. Summing over all points of $S$ gives~\eqref{eq:hodge-degree}.
\end{proof}

\subsection{Construction of the principally polarized abelian scheme $Z$}

We would like to apply Shepherd-Barron's results~\cite{sb1,sb2} to the period curve of $A$ to bound $\kappa$, but the given polarization on $A$ need not be principal. We therefore use Zarhin's construction~\cite[Theorem~1.1]{zarhin} to pass to a principally polarized abelian scheme $Z$ isogenous to $A^8$.

For a polarization $E$, let $\exp(\ker E)$ denote the least positive
integer that annihilates its kernel.
\begin{lemma}\label{lem:dyad}
Over $L^\circ$, there exists an abelian scheme $A'$ isogenous to $A$ equipped with a primitive polarization $E$ such that $4 \nmid \exp(\ker E)$.
\end{lemma}
\begin{proof}
Write $\Lambda = H_1(A,\Z)$ for the integral homology local system over $L^\circ$, and scale its alternating Riemann form $E_0 \colon \Lambda \times \Lambda \to \Z$ by its greatest common divisor so that $E_0$ is primitive. Inside $\Lambda_\Q = \Lambda \otimes_\Z \Q$, consider the dual lattice
\begin{equation*}
  \Lambda^\# = \bigl\{x \in \Lambda_\Q : E_0(x, \Lambda) \subset \Z\bigr\},
\end{equation*}
so that $\Lambda^\# / \Lambda \simeq \ker E_0$. Put $a = v_2\bigl(\exp(\Lambda^\# / \Lambda)\bigr)$. If $a \le 1$, we take $A' = A$ and $E = E_0$.

Suppose $a \ge 2$, and set
\begin{equation*}
  m = \lfloor a/2 \rfloor \ge 1, \qquad \Lambda' = 2^{-m}\Lambda \cap \Lambda^\#, \qquad E_1 = 2^m E_0.
\end{equation*}
The quotient $\Lambda' / \Lambda = (\Lambda^\# / \Lambda)[2^m]$ is the canonical $2^m$-torsion subsheaf of the polarization kernel, so it defines an $L^\circ$-isogeny from $A$ to an abelian scheme $A_1$ with homology lattice $\Lambda'$. Choose a basis $e_i,f_i$ of $\Lambda\otimes\Z_2$ such that
\[
E_0(e_i,f_j)=2^{d_i}\delta_{ij},
\qquad
E_0(e_i,e_j)=E_0(f_i,f_j)=0,
\qquad 0\le d_i\le a.
\]
Then $\Lambda^\#\otimes\Z_2$ has basis
$2^{-d_i}e_i,2^{-d_i}f_i$. Hence
$\Lambda'\otimes\Z_2=(2^{-m}\Lambda\otimes\Z_2)
\cap(\Lambda^\#\otimes\Z_2)$ has basis
$e_i'=2^{-r_i}e_i,f_i'=2^{-r_i}f_i$, where
$r_i=\min(m,d_i)$. Since $E_1=2^mE_0$, we obtain
\[
E_1(e_i',f_j')=2^{d_i+m-2r_i}\delta_{ij},
\qquad
E_1(e_i',e_j')=E_1(f_i',f_j')=0.
\]
Thus the new elementary divisors have $2$-adic valuations
\[
d_i+m-2r_i=|d_i-m|
\le \lceil a/2\rceil<a.
\]
At odd primes, the lattice is unchanged and $2^m$ is a unit. Thus $E_1$ is an integral polarization on $A_1$. Dividing $E_1$ by its integer content to make it primitive can only decrease the largest $2$-adic valuation. Iterating this step finitely many times terminates at $a \le 1$, as required.
\end{proof}

Let $(A', E)$ be the primitively polarized abelian scheme over $L^\circ$ provided by Lemma~\ref{lem:dyad}. Let $e = \exp(\ker E)$ be the exponent of its polarization kernel, with $e = 1$ if $\ker E = 0$. Set
\begin{equation*}
  r =
  \begin{cases}
    e - 1, & \text{if } e > 1,\\
    2, & \text{if } e = 1,
  \end{cases}
\end{equation*}
so that $e \mid (r+1)$ in both cases. There exists an integral quaternion matrix $J \in M_4(\Z)$ satisfying
\begin{equation*}
  J J^{\mathsf{t}} = J^{\mathsf{t}} J = r I_4.
\end{equation*}

Following Zarhin~\cite{zarhin}, consider the $L^\circ$-isogeny
\begin{equation}\label{eq:zarhin-isogeny}
  \Phi \colon (A')^4 \times (A')^4 \longrightarrow Z \coloneqq \bigl((A')^\vee\bigr)^4 \times (A')^4, \qquad (x,y) \longmapsto (Ex, \, Jx - y).
\end{equation}
The kernel of $\Phi$ is the graph of $J$ on $(\ker E)^4$. The identity $J^{\mathsf{t}} J = r I_4$ and the divisibility $e \mid (1+r)$ make $\ker \Phi$ totally isotropic with respect to the product polarization $E^{\oplus 4} \oplus E^{\oplus 4}$ with order $|\ker E|^4$. Consequently, the product polarization descends along $\Phi$ to a principal polarization on $Z$~\cite{zarhin}. Since $Z$ is isogenous to $A^8$, its normalized Hodge degree is $8h$.

\subsection{Bounding the ramification via the period map of $Z$}

Let $\overline{\calA}_{8n}^{\mathrm{perf}}$ denote the perfect cone (first Voronoi) toroidal compactification of the moduli space $\calA_{8n}$ of principally polarized abelian varieties of dimension $8n$, with boundary divisor $D = \overline{\calA}_{8n}^{\mathrm{perf}} \setminus \calA_{8n}$.

Let $p\colon L\to\overline{\calA}_{8n}^{\mathrm{perf}}$ be the extended period map of $Z$, and let $\lambda$ be the first Chern class of the determinant of the Hodge bundle. We compute the local multiplicities of $p^*D$ and prove that $\deg_L p^*(12\lambda-D)>0$. Together these give the upper bound for $\kappa$.

We first determine the local monodromy at a type $\cI$ point.

\begin{lemma}\label{lem:translation}
At a type $\cI$ point $s$, put $N=T_s-\operatorname{Id}$ on $\Lambda=H_1(A',\Z)$. Let $v$ be a primitive generator of $N(\Lambda_\Q)\cap\Lambda$, and define $k\ge1$ by
\[
  N(\Lambda)=k\Z v.
\] 
There is an integral symplectic basis
$a_1,\ldots,a_{8n},b_1,\ldots,b_{8n}$ of $H_1(Z_t,\Z)$ in which
the local monodromy fixes all $a_j$ and acts on the $b_j$ by
\[
  b_j\longmapsto
  \begin{cases}
    b_j+k a_j, & 1\le j\le 8,\\
    b_j,       & 9\le j\le 8n.
  \end{cases}
\]
\end{lemma}
\begin{proof}
Put $\delta=\gcd E(v,\Lambda)$. Since $v$ is primitive, the class of $v/\delta\in\Lambda^\#$ in $\Lambda^\#/\Lambda$ has order $\delta$, and hence $\delta\mid e$. Let $T^\vee$ denote the local monodromy on $H_1((A')^\vee,\Q)=\Lambda^\vee_\Q$.
The polarization identifies $\operatorname{Im}(T^\vee-\operatorname{Id})$ with $\Q E(v,-)$, and the definition of $\delta$ gives 
\[
\operatorname{Im}(T^\vee-\operatorname{Id})\cap\Lambda^\vee
=\Z\cdot\frac{E(v,-)}{\delta}.
\] 
Let $H=H_1((A')^8,\Z)$ and $H_Z=H_1(Z,\Z)$, and denote their local monodromy operators by $T'$ and $T_Z$. Our goal is to compute $T_Z-\operatorname{Id}$. Set 
\[
V'=\operatorname{Im}(T'-\operatorname{Id})\cap H,
\qquad
V_Z=\operatorname{Im}(T_Z-\operatorname{Id})\cap H_Z,
\] 
where the images are taken in rational homology. Since $\Phi_*$ commutes with monodromy, it maps $V'$ into $V_Z$. The first lattice has a basis consisting of $v$ in each of the eight factors. The second has a basis consisting of $w=E(v,-)/\delta$ in the first four factors and $v$ in the last four. In these bases, 
\[
\Phi_*(av,bv)=(\delta a\,w,(Ja-b)v),
\qquad a,b\in\Z^4,
\] 
so the map $\Phi_*\colon V'\to V_Z$
is represented by 
\[
P=
\begin{pmatrix}
\delta I_4&0\\
J&-I_4
\end{pmatrix}.
\] 

Choose $u\in\Lambda$ with $E(v,u)=\delta$. Then $[u]$ generates $\Lambda/v^\perp$, where $v^\perp=\{x\in\Lambda:E(v,x)=0\}$.
For the product polarization on $(A')^8$, the pairing between $(\Z v)^{\oplus8}$ and $(\Lambda/v^\perp)^{\oplus8}$ therefore
has matrix $\delta I_8$ in the bases induced by $v$ and $[u]$.

Since the polarization on $Z$ is principal and $V_Z$ is saturated, the pairing identifies $H_Z/V_Z^\perp$ with $\operatorname{Hom}(V_Z,\Z)$. Choosing the basis of this quotient dual to the preceding basis of $V_Z$, the pairing matrix is $I_8$. Let $Q$ be the matrix of the map induced by $\Phi_*$ on $H/V^{'\perp}\to H_Z/V_Z^\perp$ in the chosen bases. Since $\Phi$ pulls the polarization on $Z$ back to the product polarization on $(A')^8$, preservation of the pairings gives  $P^{\mathsf t}Q=\delta I_8$. Hence
\[
Q=\delta P^{-\mathsf t}
=
\begin{pmatrix}
I_4&J^{\mathsf t}\\
0&-\delta I_4
\end{pmatrix}.
\]
Since $\ker(T_s-\operatorname{Id})=v^\perp$ and $\operatorname{Im}(T_s-\operatorname{Id})=k\Z v$ on integral homology, we have $(T_s-\operatorname{Id})u=cv$ for some $c=\pm k$. Thus the induced map $H/V{'^\perp}\to V'$ has matrix $cI_8$. Therefore, the map $H_Z/V_Z^\perp\to V_Z$ induced by $T_Z-\operatorname{Id}$ has matrix
\[
P(cI_8)Q^{-1}=\frac{c}{\delta}PP^{\mathsf t}.
\]
Dividing by $c$, we obtain
\begin{equation}\label{eq:cusp-matrix}
  \frac{1}{\delta}PP^{\mathsf t}
  =
  \begin{pmatrix}
    \delta I_4 & J^{\mathsf t}\\
    J & \dfrac{r+1}{\delta}I_4
  \end{pmatrix}.
\end{equation}
This matrix is positive definite. It is integral since $\delta\mid e$
and $e\mid(r+1)$, and it has determinant one since $|\det P|=\delta^4$.
If $e=1$, its first diagonal entry is one. If $e>1$, then $r+1=e$,
and at least one of $\delta$ and $e/\delta$ is odd since
$4\nmid e$ by Lemma~\ref{lem:dyad}. Thus the matrix defines an odd
positive definite unimodular lattice of rank eight, which is isometric
to the standard lattice $I_8$ by~\cite{griess}.

Choosing the orientation of the monodromy $T_Z$ appropriately and using $\operatorname{ker}(T_Z-\operatorname{Id})=V_Z^{\perp}$, we obtain the conclusion.

\end{proof}

This formula determines the vanishing orders of the boundary coordinates along the period map.

\begin{lemma}\label{lem:boundary}
At a type $\cI$ point $s$, we have $\operatorname{mult}_s(p^*D)=8k$, with $k$ as in Lemma~\ref{lem:translation}. At a point of type $\cII$, $\cIII$, or $\cIV$, this multiplicity is zero.
\end{lemma}
\begin{proof}
At a point of type $\cII, \cIII$, or $\cIV$, $Z$ has good reduction after a finite base change. The extended period map on the covering curve therefore avoids $D$ at every point above $s$, so the normalized boundary contribution at $s$ is zero.

Suppose now that $s$ has type $\cI$. To avoid computing the multiplicity between period of $L$ and orbifold points, we need level structure to force all automorphisms trivial. We choose a finite cover $\pi:\widetilde L\to L$, which gives $Z$ full level $m$ structure over $\pi^{-1}(L^\circ)$, with $m\ge3$.
The extended period map of the pulled-back family and the finite map
forgetting the level structure are
\[
  \widetilde L
  \xrightarrow{\ \widetilde p\ }
  \overline{\calA}_{8n}^{\mathrm{perf}}(m)
  \xrightarrow{\ \rho\ }
  \overline{\calA}_{8n}^{\mathrm{perf}},
  \qquad \rho\circ\widetilde p=p\circ\pi.
\]
Write $D_m$ for the reduced boundary at level $m$. Since $\rho$ has ramification index $m$ along each boundary divisor, $\rho^*D=mD_m$.

Let $q_1,\ldots,q_8$ be the local 
boundary coordinates corresponding
to the rays generated by $a_1\otimes a_1,\ldots,a_8\otimes a_8$, respectively. The reduced boundary has local equation $q_1\cdots q_8=0$.

Fix $\tilde s\in\pi^{-1}(s)$ and put $d_i=\operatorname{ord}_{\tilde s}(q_i\circ\widetilde p)$.

At level $m$, a positive loop around the divisor $q_i=0$ sends $b_i$ to $b_i+m a_i$, since $m$ is the least positive $j$ for which the translation $b_i\mapsto b_i+j a_i$ acts trivially modulo $m$. Loops around the other seven divisors fix $b_i$. Thus the monodromy around $\tilde s$ sends $b_i$ to $b_i+m d_i a_i$.

Let $a$ be the ramification index of $\pi$ at $\tilde s$. The same loop projects to a loop winding $a$ times around $s$, so Lemma~\ref{lem:translation} gives $b_i\mapsto b_i+ak a_i$. Comparing the two expressions gives $m d_i=ak$ for $1\le i\le8$.

Since $\pi^*(p^*D)=\widetilde p^*(\rho^*D)$ and $\rho^*D=mD_m$, combining that the multiplicities at $\widetilde{s}$ over $s$ is $a$, we obtain  
\[
  \operatorname{mult}_s(p^*D)
  =\frac{m}{a}\sum_{i=1}^8 d_i
  =8k.
\]
\end{proof}

\begin{proposition}\label{prop:strict}
For $n\ge5$, the extended period curve of $Z$ over a very general line $L\subset\PP^n$ has strictly positive intersection number with $12\lambda-D$ on $\overline{\calA}_{8n}^{\mathrm{perf}}$.
\end{proposition}
\begin{proof}
By Shepherd-Barron's results~\cite{sb1,sb2}, the divisor $12\lambda-D$ is nef, so $\deg_L p^*(12\lambda-D)\ge0$. If $p$ is nonconstant and this degree is zero, its image curve lies in the exceptional locus of $12\lambda-D$, whose intersection with $\calA_{8n}$ is the locus of principally polarized abelian varieties with an elliptic product factor. It therefore suffices to show that $f$ is non-isotrivial and that a very general smooth fibre is simple. Indeed, the first assertion makes the period map of $A$ nonconstant on a very general line, and hence that of $Z$ as well, since $\mathbb V_Z\simeq\mathbb V_A^{\oplus8}$. The second assertion, together with $Z_t\sim A_t^8$ and $\dim A_t=n>1$, excludes elliptic isogeny factors from a very general $Z_t$. Thus the image curve meets the complement of the exceptional locus.

We first show that $f$ is non-isotrivial. More precisely, the moduli map of the smooth fibres of $f\colon X\to\PP^n$ is non-constant. If $f$ is isotrivial, its monodromy preserves a fixed polarized Hodge structure and hence has finite image, which rules out points of type $\cI$. Since $\PP^n$ is simply connected, $\pi_1(U)$ is normally generated by the meridians around the divisorial components of $\Delta$. The singular fibres over general points of these components have types $\cII$, $\cIII$, or $\cIV$. The corresponding local monodromy operators act on $H^{1,0}$ of a smooth fibre by complex reflections of orders $6$, $4$, or $3$, respectively. By~\cite[Proposition~5.4]{varmon}, this $n$-dimensional complex representation of $\pi_1(U)$ is irreducible. By the Shephard--Todd classification~\cite{st}, every primitive finite irreducible complex reflection group of rank at least five has only reflections of order two. Thus the monodromy group would have to belong to the family $G(m,p,n)$, the imprimitive finite complex reflection groups. In their standard monomial representations, every non-diagonal reflection has order two. Hence all the generating reflections, whose orders are $3$, $4$, or $6$, would be diagonal. The group they generate would then preserve each coordinate line, a contradiction to the irreducibility. Thus $f$ is non-isotrivial.

Suppose that a very general smooth fibre $X_t$ is not simple. By~\cite[Corollary~1.3 and the beginning of \S~5]{varmt},
\[
H^1(X_t,\Q)=U_1\oplus\cdots\oplus U_k,\qquad k>1,
\]
where the $U_i$ are simple Hodge structures permuted by monodromy. The irreducibility of the monodromy representation makes this
permutation action transitive. Since conjugates of divisorial meridians generate the monodromy group, some local monodromy
operator $T$ permutes the summands nontrivially.
If $T$ cyclically permutes $\ell\ge2$ summands of dimension $d$, its fixed subspace on their sum has dimension at most $d$, so
\[
(\ell-1)d\le\rk(T-\operatorname{Id})\le2.
\]
Since $d$ is positive and even, $\ell=d=2$.
On these two summands, write $T(x,y)=(Ay,Bx)$. The rank bound forces the fixed subspace to have dimension two, hence $AB=\operatorname{Id}$. Thus $T$ has eigenvalue $-1$, a contradiction to the local monodromy descriptions for types
$\cI$, $\cII$, $\cIII$, and $\cIV$. Therefore a very general fibre is simple.
\end{proof}

\begin{proof}[Proof of Theorem~\ref{thm:main}]
Let $S_{\cI} \subset S$ denote the set of type $\cI$ points. For each $s \in S_{\cI}$, write $k_s \ge 1$ for the integral monodromy width on $H_1(A',\Z)$ as in Lemma~\ref{lem:boundary}. Combining Lemma~\ref{lem:boundary} and Proposition~\ref{prop:strict} with the normalized Hodge degree $8h$ of $Z$ gives
\begin{equation*}
  8 \sum_{s \in S_{\cI}} k_s < 12 \cdot 8h.
\end{equation*}
Dividing by $8$ and substituting the formula for $h$ from Lemma~\ref{lem:degree}, we obtain the strict integer inequality
\begin{equation}\label{eq:integer-upper}
  \sum_{s \in S_{\cI}} k_s + 2 N_{\cII} + 3 N_{\cIII} + 4 N_{\cIV} < 12(n + 1).
\end{equation}
Since $k_s \ge 1$ for each $s \in S_{\cI}$, the left-hand side of~\eqref{eq:integer-upper} is at least the ramification $\kappa = N_{\cI} + 2(N_{\cII} + N_{\cIII} + N_{\cIV})$. Because both sides of~\eqref{eq:integer-upper} are integers, it follows that
\begin{equation*}
  \kappa \le 12n + 11.
\end{equation*}
Combining this upper bound with Proposition~\ref{prop:lower} yields
\begin{equation*}
  b_2(X) + 4n - 2 \le \kappa \le 12n + 11,
\end{equation*}
and hence $b_2(X) \le 8n + 13$.
\end{proof}

\renewcommand{\bibliofont}{\small}

\end{document}